\documentclass[12pt]{article}
\usepackage{amsmath}
\usepackage{amssymb}
\usepackage{amsthm}
\usepackage{tikz}
\usepackage{pgflibraryshapes}
\usepackage{url}
\newtheorem{theorem}{Theorem}

\newtheorem{conjecture}[theorem]{Conjecture}

\newtheorem{lemma}[theorem]{Lemma}
\newtheorem{question}[theorem]{Question}

\def\sqr#1#2{{\vbox{\hrule height.#2pt
    \hbox{\vrule width.#2pt height#1pt \kern#1pt
        \vrule width.#2pt}\hrule height.#2pt}}}
\def\eqed{\sqr53}
\def\qed{%
    \ifmmode\eqno\eqed
    \else\nobreak\ \hfill\eqed\medbreak\fi}

\title{Graphs with chromatic roots in the interval $(1,2)$}
\author{Gordon F. Royle\\
School of Computer Science \& Software Engineering\\
University of Western Australia\\
Nedlands WA 6009, Australia\\
{\tt gordon@csse.uwa.edu.au}
}

\begin{document}
\maketitle

\begin{abstract}
We present an infinite family of 3-connected non-bipartite graphs with chromatic roots
in the interval $(1,2)$ thus resolving a conjecture of Jackson's in the negative. In addition,
we briefly consider other graph classes that are conjectured to have 
no chromatic roots in $(1,2)$.
\end{abstract}

\section{Introduction}

The {\em chromatic polynomial} of a graph $G$ is the function $P(G,k)$ that counts
the number of $k$-colourings of $G$ when $k$ is a natural number. It is well known that $P(G,k)$ is a polynomial with degree equal to the number of vertices of $G$, and there is an extensive literature on the relationship between the properties of a graph and its chromatic polynomial. For general backgound information and references to much of this literature we refer the reader to the recent book on chromatic polynomials by Dong, Koh \& Teo \cite{MR2159409}.

One aspect of the study of chromatic polynomials that has received considerable recent interest is the problem of the distribution of the {\em chromatic roots} of a graph --- that is, the real and complex zeros of the chromatic polynomial. In a comprehensive survey paper, Jackson \cite{MR2005532} describes much  of this work and presents a number of open problems and conjectures. One of these conjectures arises from attempts to extend the portion of the real line on which the behaviour of the chromatic polynomial is completely understood. In particular, we have complete knowledge of the real roots of $P(G,x)$ for all values of $x \leq 32/27$, as described in the following fundamental theorem.

\begin{theorem}[See Jackson \cite{MR2005532}]
If $G$ is a loopless graph with $n$ vertices, $c$ components and $b$ blocks which are not isolated vertices, then 
\renewcommand{\labelenumi}{(\alph{enumi})}
\begin{enumerate}
\setlength{\itemsep}{0pt}
\item $P(G,x)$ is non-zero with sign $(-1)^n$ for $x \in (-\infty, 0)$;
\item $P(G,x)$ has a zero of multiplicity $c$ at $x=0$;
\item $P(G,x)$ is non-zero with sign $(-1)^{n+c}$ for $x \in (0,1)$;
\item $P(G,x)$ has a zero of multiplicity $b$ at $x=1$;
\item $P(G,x)$ is non-zero with sign $(-1)^{n+c+b}$ for $x \in (1,32/27] \approx (1,1.185]$.\qed
\end{enumerate}
\end{theorem}

For general graphs, the constant $32/27$ cannot be replaced by anything larger because Jackson \cite{MR1264037} produced a family of graphs with real chromatic roots arbitrarily close to $32/27$.
% and
%subsequently Thomassen \cite{MR1483433} proved that real chromatic roots are dense in 
%$[32/27,\infty)$. 
However for certain classes of graphs the chromatic-root-free interval can 
be extended --- for example, Thomassen \cite{MR1794692} showed that graphs with a 
hamiltonian path have no chromatic roots in $(1, 1.29559\ldots)$. Both Jackson's graphs
with chromatic roots close to $32/27$ and Thomassen's graphs with chromatic roots close
to $1.29559$ have cutsets of size two, raising the question of whether the chromatic-root-free interval
can be significantly extended for 3-connected graphs.
Elementary continuity arguments show that bipartite graphs of odd order have a chromatic 
root in $(1,2)$ and these can have arbitrarily high connectivity, but Jackson conjectured that
these were the only such examples.

\begin{conjecture}[Jackson \cite{MR1264037, MR2005532}]
A 3-connected graph that is not bipartite of odd order has no chromatic roots in $(1,2)$. \qed
\end{conjecture}

In this paper we demonstrate that this conjecture is false, by providing an infinite family of 3-connected non-bipartite graphs with chromatic roots in $(1,2)$, and then briefly discuss some alternative conjectures.

\section{An infinite family of counterexamples}

\begin{figure}
\begin{center}
\begin{tikzpicture}[scale=0.8]
\tikzstyle{vertex}=[circle, draw=black,inner sep = 1mm]
\node[vertex, label=270:$v_3$] (v0) at (2,1) {};
\node[vertex, label=90:$v_1$] (v1) at (2,7) {};
\node[vertex, label=270:$v_4$] (v2) at (4,1) {};
\node[vertex, label=90:$v_2$] (v3) at (4,7) {};

\node[vertex] (v4) at (-1,3) {};
\node[vertex] (v5) at (-1,4) {};
\node[vertex] (v6) at (-1,5) {};

\node[vertex] (v7) at (7,3) {};
\node[vertex] (v8) at (7,4) {};
\node[vertex] (v9) at (7,5) {};

\node[vertex, label=90:$v_0$] (v10) at (3,4) {};

\draw (v0) -- (v2);
\draw (v1) -- (v3);

\draw (v0) -- (v4);
\draw (v0) -- (v5);
\draw (v0) -- (v6);

\draw (v1) -- (v4);
\draw (v1) -- (v5);
\draw (v1) -- (v6);

\draw (v2) -- (v7);
\draw (v2) -- (v8);
\draw (v2) -- (v9);

\draw (v3) -- (v7);
\draw (v3) -- (v8);
\draw (v3) -- (v9);

\draw (v10) -- (v4);
\draw (v10) -- (v5);
\draw (v10) -- (v6);

\draw (v10) -- (v7);
\draw (v10) -- (v8);
\draw (v10) -- (v9);

\draw [dashed] (v5) ellipse (1cm and 2cm);
\draw [dashed] (v8) ellipse (1cm and 2cm);

\node at (-1,6.5) {$S$};
\node at (7,6.5) {$T$};

\end{tikzpicture}
\caption{The graph $X(3,3)$}
\label{fig1}
\end{center}
\end{figure}
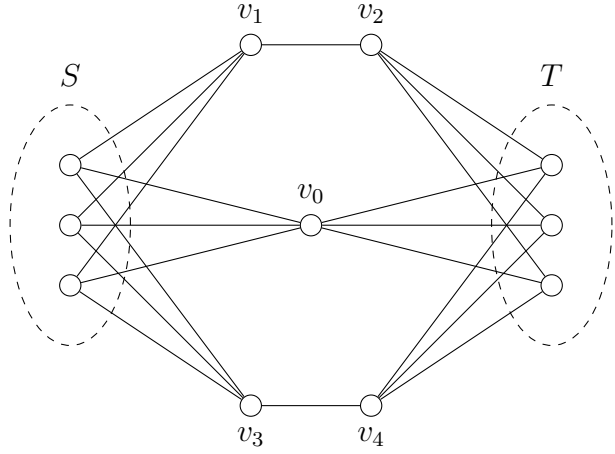

We consider a family of graphs whose smallest member is the 11-vertex graph $X(3,3)$ shown in 
Figure~\ref{fig1}. It is clear by inspection that $X(3,3)$ is 3-connected and not bipartite and it is
sufficiently small that its chromatic polynomial can be easily computed. The result of this
computation is that $P(X(3,3),x) = x(x-1)(x-2)Q(x)$, where
\begin{equation*}
Q(x) = \left( {x}^{8}-17\,{x}^{7}+
137\,{x}^{6}-677\,{x}^{5}+2228\,{x}^{4}-4969\,{x}^{3}+7284\,{x}^{2}-
6363\,x+2509 \right),
\end{equation*}
which has real roots at $r_1 \approx 1.90263148$ and $r_2 \approx 2.42196189$. Therefore 
we conclude that Jackson's conjecture is false.

Now let $X(s,t)$ be the graph obtained by replacing the two independent sets of size three labelled $S$ and $T$ in Figure~\ref{fig1} by independent sets of size $s\geq 3$ and $t \geq 3$ respectively. Our aim is to prove the following theorem:

\begin{theorem}
The graph $X(s,t)$ is 3-connected and not bipartite and if $s$, $t \geq 3$ are both odd, then it has
a chromatic root in $(1,2)$. \qed
\end{theorem}

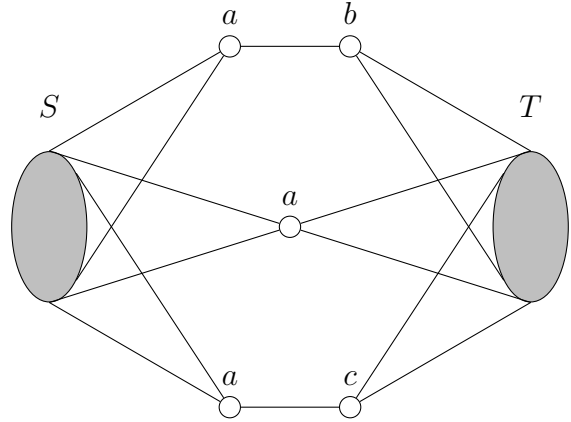
\begin{figure}
\begin{center}
\begin{tikzpicture}[scale=0.8]
\tikzstyle{vertex}=[circle, draw=black,inner sep = 1mm]
\node[vertex,  label=90:$a$] (v0) at (2,1) {};
\node[vertex,  label=90:$a$] (v1) at (2,7) {};
\node[vertex,  label=90:$c$] (v2) at (4,1) {};
\node[vertex, label=90:$b$] (v3) at (4,7) {};

%\node[vertex] (v4) at (-1,3) {};
%\node[vertex] (v5) at (-1,4) {};
%\node[vertex] (v6) at (-1,5) {};

%\node[vertex] (v7) at (7,3) {};
%\node[vertex] (v8) at (7,4) {};
%\node[vertex] (v9) at (7,5) {};

\node[vertex,label=90:$a$] (v10) at (3,4) {};

\draw (v0) -- (v2);
\draw (v1) -- (v3);

\node [shape=ellipse, minimum height = 2cm, minimum width = 1cm, fill=lightgray,draw] (A) at (-1,4) {};
\node [shape=ellipse, minimum height = 2cm, minimum width = 1cm, fill=lightgray,draw] (B) at (7,4) {};

\draw (v0) -- (A.north east);
\draw (v0) -- (A.south);

\draw (v1) -- (A.north);
\draw (v1) -- (A.south east);

\draw (v2) -- (B.north west);
\draw (v2) -- (B.south);

\draw (v3) -- (B.south west);
\draw (v3) -- (B.north);

\draw (v10)-- (A.north);
\draw (v10)-- (A.south);
\draw (v10)-- (B.north);
\draw (v10)-- (B.south);

%\node [outer sep = 0.2cm] at (v0.north) {$a$};

%\draw (v1) -- (v4);
%\draw (v1) -- (v5);
%\draw (v1) -- (v6);

%\draw (v2) -- (v7);
%\draw (v2) -- (v8);
%\draw (v2) -- (v9);

%\draw (v3) -- (v7);
%\draw (v3) -- (v8);
%\draw (v3) -- (v9);

%\draw (v10) -- (v4);
%\draw (v10) -- (v5);
%\draw (v10) -- (v6);

%\draw (v10) -- (v7);
%\draw (v10) -- (v8);
%\draw (v10) -- (v9);

\node at (-1,6) {$S$};
\node at (7,6) {$T$};

\end{tikzpicture}
\caption{One type of colouring of $X(s,t)$}
\label{fig2}
\end{center}
\end{figure}

It is clear that $X(s,t)$ is 3-connected and not bipartite, and so all that remains is to prove the claim about its chromatic roots. The chromatic polynomial of $X(s,t)$ can be calculated by considering all the possible {\em types} of colouring of the five vertices $\{v_0, v_1, v_2, v_3, v_4\}$ outside $S$ and $T$. For example, Figure~\ref{fig2} represents colourings where three distinct colours $a$, $b$ and $c$ are used,  $v_0$, $v_1$ and $v_3$ are coloured $a$ and $v_3$ and $v_4$ are coloured $b$ and $c$ respectively.
If we have a palette of $x$ colours to choose from, then there are $x$ choices for $a$, $x-1$ choices for $b$ and $x-2$ choices for $c$. Then the vertices of $S$ can be coloured anything other than $a$ and so there are $(x-1)^s$ possibilities and the vertices of $T$ can be coloured anything other than $\{a,b,c\}$ leaving $(x-3)^t$ possibilities. Thus the colourings of this particular type contribute
\begin{equation}\label{term}
x(x-1)(x-2) (x-1)^s (x-3)^t
\end{equation}
to the chromatic polynomial.

It is straightforward, though tedious, to confirm that there are exactly 27 different types of colouring for these five vertices, and therefore the chromatic polynomial of $X(s,t)$ is the sum of 27 terms of 
a similar form to the one shown in (\ref{term}), although some of terms can be combined due
to the symmetries of the situation. Overall however, the final expression is both complicated and unenlightening and more suitable for a symbolic algebra package such as Maple than for presentation in this paper. Fortunately however, we only need a small portion of the chromatic polynomial in order to deduce enough information for our purposes.

\begin{lemma}\label{derivative}
The derivative of $P(X(s,t), x)$ evaluated at $x=2$ is 
$$
P'(X(s,t), 2) = 2 \left( (-1)^s + (-1)^t + (-1)^{s+t} \right).
$$
In particular, when $s$ and $t$ are both odd
$$
P'(X(s,t), 2) = -2.
$$
\end{lemma}

\begin{figure}
\begin{center}
\begin{tikzpicture}[scale=0.8]
\tikzstyle{vertex}=[circle, draw=black,inner sep = 1mm]
\node[vertex,  label=90:$b$] (v0) at (2,1) {};
\node[vertex,  label=90:$c$] (v1) at (2,7) {};
\node[vertex,  label=90:$c$] (v2) at (4,1) {};
\node[vertex, label=90:$b$] (v3) at (4,7) {};

%\node[vertex] (v4) at (-1,3) {};
%\node[vertex] (v5) at (-1,4) {};
%\node[vertex] (v6) at (-1,5) {};

%\node[vertex] (v7) at (7,3) {};
%\node[vertex] (v8) at (7,4) {};
%\node[vertex] (v9) at (7,5) {};

\node[vertex,label=90:$a$] (v10) at (3,4) {};

\draw (v0) -- (v2);
\draw (v1) -- (v3);

\node [shape=ellipse, minimum height = 2cm, minimum width = 1cm, fill=lightgray,draw] (A) at (-1,4) {};
\node [shape=ellipse, minimum height = 2cm, minimum width = 1cm, fill=lightgray,draw] (B) at (7,4) {};

\draw (v0) -- (A.north east);
\draw (v0) -- (A.south);

\draw (v1) -- (A.north);
\draw (v1) -- (A.south east);

\draw (v2) -- (B.north west);
\draw (v2) -- (B.south);

\draw (v3) -- (B.south west);
\draw (v3) -- (B.north);

\draw (v10)-- (A.north);
\draw (v10)-- (A.south);
\draw (v10)-- (B.north);
\draw (v10)-- (B.south);

%\node [outer sep = 0.2cm] at (v0.north) {$a$};

%\draw (v1) -- (v4);
%\draw (v1) -- (v5);
%\draw (v1) -- (v6);

%\draw (v2) -- (v7);
%\draw (v2) -- (v8);
%\draw (v2) -- (v9);

%\draw (v3) -- (v7);
%\draw (v3) -- (v8);
%\draw (v3) -- (v9);

%\draw (v10) -- (v4);
%\draw (v10) -- (v5);
%\draw (v10) -- (v6);

%\draw (v10) -- (v7);
%\draw (v10) -- (v8);
%\draw (v10) -- (v9);

\node at (-1,6) {$S$};
\node at (7,6) {$T$};

\end{tikzpicture}
\caption{A type of colouring that contributes to the derivative}
\label{fig3}
\end{center}
\end{figure}
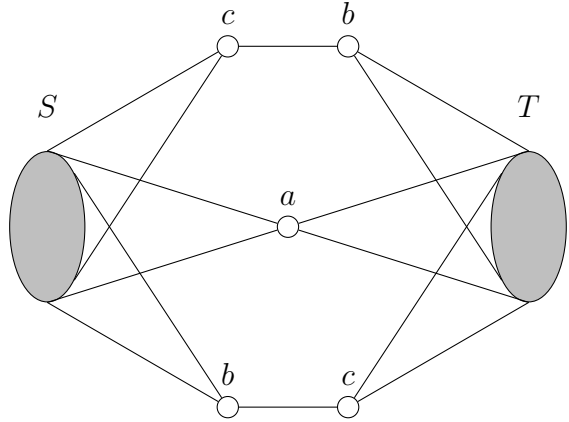

\begin{proof}
The derivative of $P(X(s,t),x)$ can be obtained by differentiating each of the 27 terms in the chromatic polynomial separately. However as our concern is only with the value of this derivative at $x=2$, we can ignore any terms with factors of $(x-2)^s$ or $(x-2)^t$ (or both) as these make no contribution to the final value (recall that $s, t > 1$). Therefore we need only consider the types of colouring that use either 1 or 3 colours on $\{v_0,v_1,v_3\}$ and similarly for $\{v_0,v_2,v_4\}$. These colourings contribute the following five terms to the chromatic polynomial:
\begin{align*}
&x \left( x-1 \right)  \left( x-2 \right)  \left( x-1 \right) ^{s}
 \left( x-3 \right) ^{t}\cr
 &x \left( x-1 \right)  \left( x-2 \right) \left( x-3 \right) ^{s} \left( x-3 \right) ^{t}\cr
 &x \left( x-1 \right)  \left( x-2 \right)  \left( x-3 \right) ^{s} \left( x-1 \right) ^{t}\cr
 &2\,x \left( x-1 \right)  \left( x-2 \right)  \left( x-3 \right) \left( x-3 \right) ^{s} \left( x-3 \right) ^{t}\cr
 &x \left( x-1 \right) \left( x-2 \right)  \left( x-3 \right)  \left( x-4 \right)  \left( x-3 \right) ^{s} \left( x-3 \right) ^{t}\cr
\end{align*}
(For example, the second term $x \left( x-1 \right)  \left( x-2 \right) \left( x-3 \right) ^{s} \left( x-3 \right) ^{t}$
counts the colourings of the type shown in Figure~\ref{fig3}.)
The sum of the last two terms is 
$$
x \left( x-1 \right)  \left( x-2 \right) ^{2} \left( x-3 \right) ^{1+s+t}
$$
and so these also make no net contribution to the derivative of $P(X(s,t),x)$ evaluated at $x=2$. Differentiating the remaining three terms and substituting $x=2$ gives the stated result, and the second
part follows immediately.
\end{proof}

When $s$ and $t$ are both odd, the graph $X(s,t)$ has an odd number of vertices and as it is 
3-connected, the chromatic root at $x=1$ has multiplicity one. Therefore $P(X(s,t),x)$ is negative
on $(-\infty, 0)$, positive on $(0,1)$ and negative on some small interval $(1,1+\alpha)$. However by Lemma~\ref{derivative}, the chromatic polynomial is {\em decreasing} as it passes through the root at $x=2$ and so is {\em positive} on some small interval $(2-\beta,2)$. It follows by continuity that the chromatic polynomial has a root in the interval $(1,2)$.

We can perform a similar analysis on the graphs obtained by adding the single edge $v_1v_4$
to $X(s,t)$ with the following result.

 \begin{theorem}
If $Y(s,t)$ is the graph obtained from $X(s,t)$ by joining $v_1$ and $v_4$ then
$$
P'(Y(s,t), 2) = 2 \left( (-1)^s + (-1)^t + (-1)^{s+t}\right)
$$
and hence when $s, t$ are both odd, the chromatic polynomial of $Y(s,t)$ also has a chromatic
root in $(1,2)$.\qed
\end{theorem}

\section{Questions and Conjectures}

The results of the previous section raise two immediate questions:
\begin{question}
What is the largest $\delta$ such that $(1,\delta)$ is a chromatic-root-free interval for 3-connected non-bipartite graphs?
\end{question}
\begin{question}
Which classes of graphs have no chromatic roots in $(1,2)$?
\end{question}
The first question is asking whether there are other 3-connected non-bipartite graphs with chromatic roots {\em smaller} than $1.90263148$ and if so, how much smaller. The first obvious place to look is to see what happens when varying $s$ and $t$ for $X(s,t)$ and $Y(s,t)$. However the proof of our main theorem is an existence proof and gives us no clues as to the actual location of the chromatic roots other than that they are in $(1,2)$.  However computational evidence (detailed in Table~\ref{tabx} and Table~\ref{taby}) suggests that $X(3,3)$ contributes the smallest chromatic roots in these families. 

A related question is obtained by {\em not} excluding bipartite graphs and just requiring 3-connectivity. In this case, the smallest {\em known} chromatic root in $(1,2)$ is approximately $1.7811$, which comes from the complete bipartite graph $K_{3,4}$.

\begin{table}
\begin{center}
{\small
\begin{tabular}{|c|ccccccccc|}
\hline
$s \backslash t$ & $3$ & $5$ & $7$ & $9$ & $11$ & $13$ & $15$ & $17$ & $19$ \\
\hline
$3$&$1.9026$&$1.9223$&$1.9342$&$1.9424$&$1.9484$&$1.9531$&$1.9568$&$1.9599$&$1.9625$\\
$5$&&$1.9372$&$1.9464$&$1.9527$&$1.9574$&$1.9611$&$1.9640$&$1.9665$&$1.9685$\\
$7$&&&$1.9539$&$1.9591$&$1.9630$&$1.9660$&$1.9685$&$1.9706$&$1.9723$\\
$9$&&&&$1.9636$&$1.9669$&$1.9696$&$1.9717$&$1.9735$&$1.9751$\\
$11$&&&&&$1.9699$&$1.9722$&$1.9741$&$1.9757$&$1.9771$\\
$13$&&&&&&$1.9743$&$1.9761$&$1.9775$&$1.9788$\\
$15$&&&&&&&$1.9777$&$1.9790$&$1.9801$\\
$17$&&&&&&&&$1.9802$&$1.9813$\\
$19$&&&&&&&&&$1.9822$\\
\hline
\end{tabular}
}
\end{center}
\caption{Chromatic root of $X(s,t)$ to 4 decimal places}
\label{tabx}
\end{table}

\begin{table}
\begin{center}
{\small
\begin{tabular}{|c|ccccccccc|}
\hline
$s \backslash t$ & $3$ & $5$ & $7$ & $9$ & $11$ & $13$ & $15$ & $17$ & $19$ \\
\hline
$3$&$1.9131$&$1.9294$&$1.9397$&$1.9468$&$1.9521$&$1.9563$&$1.9596$&$1.9624$&$1.9648$\\
$5$&&$1.9420$&$1.9500$&$1.9556$&$1.9598$&$1.9631$&$1.9659$&$1.9681$&$1.9700$\\
$7$&&&$1.9566$&$1.9613$&$1.9648$&$1.9676$&$1.9699$&$1.9718$&$1.9734$\\
$9$&&&&$1.9653$&$1.9684$&$1.9708$&$1.9728$&$1.9745$&$1.9759$\\
$11$&&&&&$1.9711$&$1.9733$&$1.9751$&$1.9766$&$1.9778$\\
$13$&&&&&&$1.9752$&$1.9768$&$1.9782$&$1.9794$\\
$15$&&&&&&&$1.9783$&$1.9796$&$1.9807$\\
$17$&&&&&&&&$1.9807$&$1.9817$\\
$19$&&&&&&&&&$1.9827$\\
\hline
\end{tabular}
}
\end{center}
\caption{Chromatic root of $Y(s,t)$ to 4 decimal places}
\label{taby}
\end{table}

With respect to the second question, a number of different classes have been considered by 
various authors. As all (previously) known examples of 2-connected graphs with chromatic roots in $(1,2)$ had cutsets of size two or were bipartite of odd order, it was natural to conjecture that these were the only such graphs. However, with the considerable advantage of hindsight, the condition ``3-connected and non-bipartite'' is somewhat unsatisfactory because an ideal condition would {\em intrinsically} exclude odd order bipartite graphs rather than explicitly omit them from consideration.

One such condition was given by Thomassen \cite{MR1426436}, who conjectured that {\em hamiltonian} graphs do not have chromatic roots in $(1,2)$. Recall that a graph $G$ is said to be {\em $1$-tough} if there is no subset $S$ of its vertices such that $G-S$ has more than $|S|$ components, and that a hamiltonian graph is necessarily $1$-tough. Now $X(s,t) $ and $Y(s,t)$ both have the property that they have sets of 3 vertices (e.g. $\{v_0,v_1,v_3\}$) whose removal leaves more than 3 components, and therefore they are not $1$-tough and hence not hamiltonian. Thus Thomassen's conjecture remains open. 

In recent work, Dong \& Koh \cite{MR2272231} have found a number of sufficient conditions for graphs to have no chromatic roots in $(1,2)$, such as the following useful result.
\begin{theorem}[Dong \& Koh \cite{MR2272231}]
If the graph $G$ has a hamiltonian path $v_1 \sim v_2 \sim \cdots \sim v_n$ such that for all $i > 2$, the vertex $v_i$ has at least one neighbour in $\{v_1,v_2, \ldots,v_{i-2}\}$ (i.e. at least two ``back neighbours'' in total including $v _{i-1}$) then $G$ has no chromatic roots in $(1,2)$. \qed
\end{theorem}
One feature common to all of the graphs that Dong \& Koh can {\em prove} to be chromatic-root-free in the interval $(1,2)$ is that they satisfy a condition that is a weaker version of $1$-toughness, namely that they have no {\em independent} sets $S$ whose removal leaves more than $|S|$ components.
%\footnote{It is very tempting to coin the phrase ``weakly 1-tough'' for this property, but I shall resist.}. 
This suggests that maybe it is the {\em toughness} of hamiltonian graphs (rather than the actual cycle itself) that is important, and along these lines Dong \& Koh have made the following interesting
strengthening of Thomassen's conjecture.

\begin{conjecture}[Dong \& Koh \cite{MR2272231}]
If $G$ is a graph with a chromatic root in $(1,2)$, then $G$ has an independent set $S$ such that
$$
c(G-S) > |S|
$$
where $c(G-S)$ denotes the number of connected components of $G-S$.
\end{conjecture}

\bibliographystyle{acm2url}
\bibliography{chromaticroots}

 \end{document}